\documentclass[11pt]{article}
\usepackage{graphicx,mathrsfs,amssymb}
\usepackage{amsmath,latexsym,amsbsy,amssymb}
\usepackage{graphicx}

\def\ds{\displaystyle}
\def\forall{\hbox{for all}~}

\def\ve{\varepsilon}

\def\R{I\!\!R}

\def\vs{\vskip 2em}

\def\v{\vskip 1em}

\def\bega{\begin{array}}
\def\enda{\end{array}}
\def\begi{\begin{itemize}}
\def\endi{\end{itemize}}

\def\bel{\begin{equation}\label}
\def\eeq{\end{equation}}
\def\sqr#1#2{\vbox{\hrule height .#2pt
\hbox{\vrule width .#2pt height #1pt \kern #1pt
\vrule width .#2pt}\hrule height .#2pt }}
\def\square{\sqr74}
\def\endproof{\hphantom{MM}\hfill\llap{$\square$}\goodbreak}

\newtheorem{theorem}{Theorem}[section]

\newtheorem{definition}[theorem]{Definition}

\newtheorem{lemma}[theorem]{Lemma}

\begin{document}
\title{\bf Covering numbers for  bounded variation functions}\vs
\author{
Prerona Dutta and  Khai T. Nguyen\\
\quad\quad
{\small Department of Mathematics}\\
{\small North Carolina State University.} \\
\quad\\
{\small Emails: pdutta@ncsu.edu,\qquad khai@math.ncsu.edu}}
\maketitle
\begin{abstract}
In this paper, we provide upper and lower estimates for the minimal number of functions needed to  represent a bounded variation function with an accuracy of epsilon with respect to ${\bf L}^1$--distance.
\end{abstract}
\section{Introduction}
\label{sec:1}
\setcounter{equation}{0}
The ${\ve}$-entropy has been studied extensively in a variety of literature and disciplines.  It plays a central role in various areas of information theory and statistics, including nonparametric function estimation, density information, empirical processes and machine learning (see e.g in \cite{LB, DH, DP}). This concept was first introduced by Kolmogorov and Tikhomirov in \cite{KT}:
\begin{definition}\label{def1}
Let $(X,d)$ be a metric space and $E$ a precompact subset of $X$. For $\varepsilon >0$, let $\mathcal{N}_{\varepsilon}(E|X)$ be the minimal number of sets in an $\varepsilon$-covering of $E$, i.e., a covering of $E$ by subsets of $X$ with diameter no greater than $2\varepsilon$. Then $\varepsilon$-entropy of $E$ is defined as 
\[
\mathcal{H}_{\varepsilon}(E~|~X)=\log_2\mathcal{N}_{\varepsilon}(E~|~X).
\]
\end{definition}
In other words, it is the minimum number of bits needed to represent a point in a given set $E$ in the space $X$ with an accuracy of $\varepsilon$ with respect to the metric $d$. 
 
 A classical topic in the field of probability is to investigate the metric covering numbers for general classes of real-valued functions $\mathcal{F}$ defined on $X$ under the family of ${\bf L}^1(dP)$ where $P$ is a probability distribution on $X$. Upper bounds in terms of  Vapnik-Chervonenkis and pseudo-dimension of the function class were established in \cite{RMD}, and then improved in \cite{DP,DH, DH1}. Several results on lower bounds were also studied in \cite{KMJ}. Later on, upper and lower estimates of the ${\ve}$-entropy of $\mathcal{F}$ in ${\bf L}^1(dP)$ in terms of a scale-sensitive dimension of the function class were provided  in \cite{LBW,KMJ}, and applied to machine learning.  
\quad\\
\quad\\
Thanks to the Helly's theorem, a set of uniformly bounded variation functions is compact in ${\bf L}^1$-space. A natural question is to quantify the compactness of such sets by using the $\ve$-entropy. In \cite{KMJ}, the authors considered this problem in the scalar case and proved  that the $\ve$-entropy of  a class of real valued functions of bounded variation in ${\bf L}^1$ is of the order of $\ds {1\over \ve}$. Some related works have been done in the context of density estimation where attention has been given to the problem of finding covering numbers for the classes of densities that are unimodal or nondecreasing  in \cite{LB, PG}. In the multi-dimensional cases, the covering numbers of convex and uniformly bounded functions  were studied in \cite{GS}. It was shown that the $\ve$-entropy of a class of convex functions with uniform bound in ${\bf L}^1$ is of the order of $\ds{1\over \ve^{n\over 2}}$ where $n$ is the dimension of the state variable. The result was previously studied for scalar state variables in \cite{DD} and for convex functions that are uniformly bounded and uniformly Lipschitz with a known Lipschitz constant in \cite{EB}. These results have direct implications in the study of rates of convergence of empirical minimization procedures (see e.g. in \cite{LB1,SVG} as well as optimal convergence rates in the numerous convexity constrained function estimation problems (see e.g. in \cite{LB0, LLC,YB}).

Recently, the ${\ve}$-entropy has been used to measure the set of solutions of certain nonlinear partial different equations. In this setting, it  could provide a measure of the order of ``resolution'' and of the ``complexity'' of a numerical scheme, as suggested in~\cite{Lax78, Lax02}. Roughly speaking, the order of magnitude of the $\varepsilon$-entropy should indicate the minimum number of operations that one should perform  in order to obtain an approximate solution with a precision of order $\varepsilon$ with respect to the considered topology. A starting point of this research topic is a result which was obtained in \cite{DLG} for  a scalar conservation law in one dimensional space
\bel{CL}
u_t(t,x)+f(u(t,x))_x~=~0,
\eeq
with uniformly convex flux $f$. They showed that the upper bound of the minimum number of functions needed to represent an entropy solution $u$ of (\ref{CL}) at any time $t>0$ with accuracy $\varepsilon$ with respect to $\bf{L}^1$-distance is of the order of  $\ds {1\over\ve }$. In \cite{AON1} a lower bound on such an $\varepsilon$-entropy was established, which is of the same order as of the upper bound in \cite{DLG}.  More generally, the authors in \cite{AON1} also obtained the same estimate for a system of hyperbolic conservation laws in \cite{AON2, AON3}.  In the scalar case, it is well-known that the integral form of an entropy solution of (\ref{CL}) is a viscosity  solution of the related Hamilton-Jacobi equation. Therefore, it is natural to study the ${\ve}$-entropy for the set of viscosity solutions to the Hamilton-Jacobi equation
\bel{HJ}
u_t(t,x)+H\big(\nabla_x u(t,x) \big)~=~0\,,
\eeq
with respect to $\bf{W}^{1,1}$-distance in multi-dimensional cases. Most recently, it has been proved in \cite{ACN} that  the minimal number of functions needed to represent a viscosity solution of (\ref{HJ})  with accuracy $\varepsilon$ with respect to the $\bf{W}^{1,1}$-distance is of the order of $\ds {1\over\varepsilon^n}$, provided that  $H$ is uniformly convex. Here, $n$ is the dimension of the state variable. The same result for when the Hamiltonian depends on the state variable $x$ has also been obtained by the same authors in \cite{ACN1}. 
\qquad\\
\quad\\
Interestingly, the authors in \cite{ACN} also  established an upper bound on the ${\ve}$-entropy for the class of monotone functions in $\bf{L}^{1}$-space. As a consequence of Poincar\'e-type inequalities, they could obtain the ${\ve}$-entropy for a class of semi-convex/concave functions in Sobolev $\bf{W}^{1,1}$ space. This result somehow extended the one in \cite{GS,DD,EB} to a stronger norm, $\bf{W}^{1,1}$-norm instead of ${\bf L}^1$-norm. Motivated by the results in \cite{KMJ,GS,DD,EB, ACN} and a possible application to Hamilton-Jacobi equation with non-strictly convex Hamiltonian, we will provide in the present paper upper and lower estimates of  the ${\ve}$-entropy for a class of uniformly bounded total variation functions in $\bf{L}^{1}$-space in multi-dimensional cases.  In particular, our result shows that the minimal number of functions needed to  represent a function with bounded variation with an error $\ve$ with respect to ${\bf L}^1$-distance is of the order of ${1\over \ve^n}$. The precise statement will be stated in Theorem \ref{main} in section 3. 

\section{Notations and preliminaries}
\setcounter{equation}{0}
Let $n\geqslant 1$ be an integer and $D$ be a measurable subset of $\R^n$. Throughout the paper we shall denote by:
\begin{itemize}
\item $|\cdot|$ the Euclidean norm in $\R^n$;
\item $\langle\cdot,\cdot\rangle$ the Euclidean inner product in $\R^n$;
\item $\mathrm{int}(D)$ the interior of $D$;
\item $\partial D$ the boundary of $D$;
\item $\mathrm{Vol}(D)$ the Lebesgue measure of a measurable set $D\subset \R^n$;
\item $\mathbf{L}^{1}(D,\R)$ the Lebesgue space of all (equivalence classes of) summable real functions on $D$, equipped with the usual norm $\|\cdot\|_{\mathbf{L}^{1}(D)}$;
\item $\mathbf{L}^{\infty}(D,\R)$  the space of all essentially bounded real functions on~$D$,  and by $\|u\|_{\mathbf{L}^{\infty}(D)}$ the essential supremum of a function $u\in \mathbf{L}^{\infty}(D,\R)$;
\item $\mathcal{C}^1_{c}(\Omega,\R^n)$, with $\Omega\subset\R^n$ an open set, the set of all continuous differentiable functions from $\Omega$ to $\R^n$ with a compact support in $\Omega$;
\item $\chi_{D}(x)=\left\{\bega{rl}
&1 \qquad~~\mathrm{if}\qquad x\in D\,,
\\[4mm]
&0\quad~\mathrm{if}\qquad x\in\R^n\backslash D\,
\enda\right.$
the characteristic function of a subset $D$ of $\R^n$.
\item $\mathrm{Card}(S)$ the number of elements of any finite set $S$;
\item $\lfloor x\rfloor\doteq a\doteq\max\{z\in\mathbb{Z}~|~z\leq x\}$ denotes the integer part of $x$.
\end{itemize}
We now introduce the concept of functions of bounded variations.
\begin{definition} The function $u\in {\bf L}^1(\Omega,\R)$ is {\it a function of bounded variation on $\Omega$ (denoted by $BV(\Omega,\R)$)}  if the distributional derivative of $u$ is representable by a finite Radon measure in $\Omega$, i.e., if  
\[
\int_{\Omega}~u\cdot {\partial \varphi\over\partial x_i}~dx~=~-\int_{\Omega}\varphi dD_iu\qquad\qquad\forall \varphi\in\mathcal{C}_c^1(\Omega,\R), i\in\{1,2,...,n\}
\]
for some Radon measure $Du=(D_1u,D_2u,...,D_nu)$. We denote by $|Du|$ the total variation of the vector measure $Du$, i.e., 
\[
|Du|(\Omega)~=~\sup\left\{\int_{\Omega}u(x)\mathrm{div}(\phi)~\Big|~\phi\in\mathcal{C}_c^1(\Omega,\R^n), \|\phi\|_{{\bf L}^{\infty}(\Omega)}\leq 1\right\}\,.
\]
\end{definition}
Let's recall a Poincar\'e-type inequality for bounded total variation functions on convex domain that will be used in the paper. This result is based on \cite[theorem 3.2]{AD} and on \cite[Proposition 3.2.1, Theorem 3.44]{ANP}.
\begin{theorem}(Poincar\'e inequality) Let $\Omega\subset \R^n$ be an open,  bounded, convex set with Lipschitz boundary. For any  $u\in BV(\Omega,\R)$, it holds
\[
\int_{\Omega} \big|u(x)-u_{\Omega}\big|~dx~\leq~{\mathrm{diam}(\Omega)\over 2}\cdot |Du|(\Omega)
\]
where 
\[
u_{\Omega}~=~{1\over \mathrm{Vol}(\Omega)}\cdot \int_{\Omega}u(x)~dx
\]
is the mean value of $u$ over $\Omega$.
\end{theorem}
To complete this section, we will state a result  on the $\ve$-entropy  for  a class of bounded total variation functions in the scalar case using a method similar to the one provided in \cite{BKP}. Given $L,V, M>0$, denote by 
\bel{B}
\mathcal{B}_{[L,M,V]}~=~\left\{f\in {\bf L}^1([0,L],[0,M])~\Big|~ |Df|((0,L))\leq V\right\}\,.
\eeq
\begin{lemma}\label{1-D-BV} For all $0<\ve<{L(M+V)\over 6}$, it holds
\bel{BV-Es1}
\mathcal{H}_{\ve}\left(\mathcal{B}_{[L,M,V]}~\big|~{\bf L}^1([0,L])\right)~\leq~8\cdot \left[{L(M+V)\over \ve}\right]\,.
\eeq
\end{lemma}
{\bf Proof.} For any $f\in \mathcal{B}_{[L,M,V]}$, let $V_f(x)$ be the total variation of $f$ over $[0,x]$. We decompose 
\[
f(x)~=~f^+(x)-f^-(x)\qquad\forall x\in [0,L]\,.
\]
where $f^{-}={V_f-f\over 2} +{M\over 2}$ is a non-decreasing function $[0,L]$ to $\left[0,{L+M\over 2}\right]$ and  $f^{+}={V_f+f\over 2} +{M\over 2}$ is a nondecreasing function $[0,L]$ to $\left[{M\over 2},{L+2M\over 2}\right]$. Denote by 
\[
\mathcal{I}~:=~\left\{g: [0,L]\to \left[0,{V+M\over 2}\right]~\Big|~g~\mathrm{is~nondecreasing}\right\}\,,
\]
we then have  
\bel{inc1}
\mathcal{B}_{[L,M,V]}~\subseteq~\left(\mathcal{I}+{M\over 2}\right)-\mathcal{I}~\doteq~\left\{g-h~\Big|~g\in\mathcal{I}+{M\over 2}\quad\mathrm{and}\quad h\in \mathcal{I}\right\}\,.
\eeq
For any $\ve>0$, it holds 
\[
\mathcal{N}_{\ve}\left(\mathcal{B}_{[L,M,V]}~|~{\bf L}^1([0,L])\right)~\leq~\left[\mathcal{N}_{{\ve\over 2}}(\mathcal{I}~|~{\bf L}^1([0,L]))\right]^2\,.
\]
Indeed, from the definition \ref{def1}, there exists a set $\mathcal{G}_{{\ve\over 2}}$ of $\mathcal{N}_{{\ve\over 2}}(\mathcal{I}~|~{\bf L}^1([0,L]))$ subsets  of ${\bf L}^{1}([0,L])$ such that 
\[
\ds\mathcal{I}~\subseteq~\bigcup_{\mathcal{E}\in\mathcal{G}_{{\ve\over 2}}}~\mathcal{E}\quad\mathrm{and}\quad\mathrm{diam}(\mathcal{E})~=~\sup_{h_1,h_2\in\mathcal{E}}\|h_1-h_2\|_{{\bf L^1([0,L])}}~\leq~\ve\,.
\]
Thus, (\ref{inc1}) implies 
\[
\mathcal{B}_{[L,M,V]}~\subseteq~\bigcup_{(\mathcal{E}_1,\mathcal{E}_2)\in\mathcal{G}_{{\ve\over 2}}\times \mathcal{G}_{{\ve\over 2}}} \left[\left(\mathcal{E}_1+{M\over 2}\right)-\mathcal{E}_2\right]\,.
\]
For any two functions
$$
f_i~=~g_i-h_i~\in \left(\mathcal{E}_1+{M\over 2}\right)-\mathcal{E}_2\qquad\mathrm{for}~i=1,2\,,
$$
we have
\[
\bega{rl}
\|f_1-f_2\|_{{\bf L}^1([0,L])}&\leq~\|g_1-g_2\|_{{\bf L}^1([0,L])}+\|h_1-h_2\|_{{\bf L}^1([0,L])}\\[3mm]
&\leq~\ds\mathrm{diam}\left(\mathcal{E}_1+{M\over 2}\right)+\mathrm{diam}(\mathcal{E}_2)~\leq~\ve+\ve~=~2\ve
\enda
\]
and this implies that 
\[
\ds \mathrm{diam}\left[\left(\mathcal{E}_1+{M\over 2}\right)-\mathcal{E}_2\right]~\leq~2\ve\,.
\]
By the definition \ref{def1}, we have 
\[
\mathcal{N}_{\ve}\left(\mathcal{B}_{[L,M,V]}~\Big|~{\bf L}^1([0.L])\right)~\leq~\mathcal{N}^2_{{\ve\over 2}}(\mathcal{I}~|~{\bf L}^1([0,L]))\,.
\]
and thus
\bel{ess1}
\mathcal{H}_{\ve}\left(\mathcal{B}_{[L,M,V]}~\Big|~{\bf L}^1([0,L])\right)~\leq~2\cdot \mathcal{H}_{\ve\over 2}\left(\mathcal{I}~\Big|~{\bf L}^1([0,L])\right)\,.
\eeq
\quad\\
Finally, applying \cite[Lemma 3.1]{DLG} for $\mathcal{I}$, we obtain that  for $0<\ve<{L(M+V)\over 6}$, it holds 
\[
\ds\mathcal{H}_{{\ve\over 2}}\left(\mathcal{I}~\big|~{\bf L}^1([0,L])\right)~\leq~4\cdot \left\lfloor{L(M+V)\over \ve}\right\rfloor\,,
\]
and (\ref{ess1}) yields (\ref{BV-Es1}).
\endproof
\section {Estimates of  the $\ve$-entropy for a class of BV functions}
\setcounter{equation}{0}
In this section, we establish upper and lower estimates of  the ${\ve}$-entropy for a class of uniformly bounded total variation functions,
\bel{F}
\mathcal{F}_{[L,M,V]}~=~\left\{u\in {\bf L}^1([0,L]^n,\mathbb{R})~\Big|~\|u\|_{{\bf L}^{\infty}([0,L]^{n})}\leq M, |Du|((0,L)^n)\leq V \right\}\,,
\eeq
 in the ${\bf L}^{1}([0,L]^n,\R)$-space. In particular, it is shown that the minimal number of functions needed to  represent a function in  $\mathcal{F}_{[L,M,V]}$ with an error $\ve$ with respect to ${\bf L}^1$-distance is of the order of ${1\over \ve^n}$. More precisely, our main result is stated as the following.
 
 \v
 
\begin{theorem}\label{main} Given $L,M,V>0$, for every $0<\ve< {ML^n\over 8}$, it holds
\bel{m-est}
{\log_2(e)\over 8}\cdot \left\lfloor{VL\over 2^{n+2}\ve}\right\rfloor^n~\leq~\mathcal{H}_{\ve}\left(\mathcal{F}_{[L,M,V]}~\Big|~{\bf L}^1([0,L]^n)\right)~\leq~\Gamma_{[n,L,M,V]}\cdot {1\over \ve^n}
\eeq
where the constant $\Gamma_{[n,L,M,V]}$ is computed as 
\[
\Gamma_{[n,L,M,V]}~=~{8\over \sqrt{n}}\left(4\sqrt{n}LV\right)^n+\left({2^{n+7}V\over M}+8\right)\cdot \left({ML^n\over 8}\right)^n\,.
\]
\end{theorem}

{\bf Proof.} {\it (Upper estimate)}~Let's first prove the upper-estimate of $\mathcal{H}_{\ve}\left(\mathcal{F}_{[L,M,V]}~\Big|~{\bf L}^1([0,L]^n)\right)$. The proof is divided into several steps:

{\bf 1.} For any $N\in\mathbb{N}$, we divide the square $[0,L]^n$ into $N^n$ small squares $\square_{\iota}$ for $\iota=(\iota_1,\iota_2,...,\iota_n)\in \{0,1,...,N-1\}^n$ such that 
\[
\square_{\iota}~=~{\iota L\over N}+\Bigg(\left[0, {L\over N}\right]\times \left[0, {L\over N}\right]\times...\times  \left[0, {L\over N}\right]\Bigg)\qquad\mathrm{and}\qquad \bigcup_{\iota\in \{0,1,2,...,N-1\}^n}~\square_{\iota}~=~[0,L]^n\,.
\]
For any $u\in\mathcal{F}_{[L,M,V]}$, denote by 
\[
-M~\leq~u_{\iota}~=~{1\over\mathrm{Vol}(\square_{\iota})}~\int_{\square_{\iota}}u(x)~dx~\leq~M
\]
the average value of $u$ in $\square_{\iota}$ for every $\iota\in \{0,1,2,...,N-1\}^n$. Let $\tilde{u}$ be a piecewise constant function on $[0,L]^n$ such that 
\begin{equation*}
\tilde{u}(x)~=~\left\{\bega{rl}
&\ds u_{\iota}~\qquad\qquad\forall x\in \mathrm{int}\big(\square_{\iota}\big)\,,
\\[4mm]
&\ds 0\qquad~~\qquad\forall x\in \bigcup_{\iota\in\{1,2,\dots,N-1\}^n}\partial\square_{\iota}\,.
\enda\right. 
\end{equation*}
Thanks to the Poincar\'e inequality, we have
\[
\int_{\square_{\iota}}|u(x)-u_{\iota}|~dx~\leq~{\mathrm{diam}(\square{\iota})\over 2}\cdot |Du|(\mathrm{int}(\square_{\iota}))
\]
for all $\iota\in \{0,1,2,...,N-1\}^n$. Hence, the ${\bf L}^1$-distance between $u$ and $\tilde{u}$ can be estimated as follows
\begin{multline}\label{L-est1}
\|u - \tilde{u}\|_{{\bf L}_1([0,L]^n)}~=~ ~\int_{[0,L]^n}|u(x)-\tilde{u}(x)|~dx~=~\sum_{\iota\in \{0,1,2,...,N-1\}^n} \int_{\square{\iota}} |u(x) - u_{\iota}|~dx\\[2mm]
~\leq~\sum_{\iota\in \{0,1,2,...,N-1\}^n}~\Bigg( {\mathrm{diam}(\mathrm{int}(\square{\iota}))\over 2}\cdot |Du|(\mathrm{int}(\square{\iota})) \Bigg)~\leq~ \frac{L\sqrt{n}}{N}~\sum_{\iota\in \{0,1,2,...,N-1\}^n} |Du|(\square{\iota})\\[2mm]
~=~{L\sqrt{n}\over N}~|Du|((0,L)^n)~\leq~{L\sqrt{n}\over N}\cdot V\,.
\end{multline}

\v

{\bf 2.} Let $e_1,e_2,...,e_n$ be the standard basis of $\R^n$ where $e_i$ denotes the vector with a $1$ in the $i$-th coordinate and  $0$'s elsewhere. For any $\iota\in\{0,1,2,...,N-1\}^n$ and $j\in \{1,2,...,n\}$, we estimate $\left|u_{\iota+e_j}-u_{\iota}\right|\,$ in the following way:

\begin{multline}\label{qs-mono}
|u_{\iota+e_j}-u_{\iota}|~=~ \left| {1\over \mathrm{Vol}\left(\square_{\iota+e_j}\right)}~\int_{\square_{\iota+e_j}}u(x)~dx  - {1\over \mathrm{Vol}\left(\square_{\iota}\right)}~\int_{\square_{\iota}}u(x)~dx
\right|\\
~=~{1\over \mathrm{Vol}\left(\square_{\iota}\right)}\cdot \left| \int_{\square_{\iota}}~u\Big(x+{L\over N}\cdot e_j\Big)-u(x)~dx\right|~=~{1\over \mathrm{Vol}\left(\square_{\iota}\right)}\cdot \left| \int_{\square_{\iota}}\int_{0}^{{L\over N}}~Du(x+se_j)(e_j)~dsdx\right|\\
~\leq~{1\over \mathrm{Vol}\left(\square_{\iota}\right)}\cdot\int_{0}^{{L\over N}}  \left|\int_{\square_{\iota}}~Du(x+se_j)(e_j)~dx\right|ds~~\leq~\left({N\over L}\right)^{n-1}\cdot |Du|(\mathrm{int}(\square_{\iota}\cup \square_{\iota+e_j} ))\,.
\end{multline}
 Let us rearrange the index set 
 \[
 \{0,1,2,\dots,N-1\}^{n}~=~\left\{\kappa^1,\kappa^2,\dots,\kappa^{N^n}\right\}
 \]
 in the way such that for all $j\in\{1,...,N^n-1\}$, it holds
 \[
 \kappa^{j+1}~=~\kappa^{j}+e_k\qquad\mathrm{for\ some}~k\in\{1,2,...,n\}\,.
 \]
From (\ref{qs-mono}) and (\ref{F}), we have 
\begin{multline}\label{TV1}
 \sum_{j=1}^{N^n}\left|u_{\kappa^{j+1}}-u_{\kappa^j}\right|~\leq~\left({N\over L}\right)^{n-1}\cdot \sum_{j=1}^{N^n}|Du|(\mathrm{int}(\square_{ \kappa^{j}}\cup \square_{ \kappa^{j+1}}))
 \cr
 ~\leq~2\left({N\over L}\right)^{n-1}\cdot |Du|((0,L)^n)~\leq~2V\left({N\over L}\right)^{n-1}\,.
\end{multline}
 To conclude this step, we define the function $f_{u,N}: [0,LN^{n-1}]\to [-M,M]$ associated with $u$ such that 
 \[
 f_{u,N}(x)~=~u_{\kappa^{i}}\qquad\forall x\in \left[{i\cdot L\over N},{(i+1)\cdot L\over N}\right), i\in\left\{1,2,...,N^{n}-1\right\}\,.
 \]
 Recalling (\ref{TV1}), we have
  \bel{TV2}
|Df_{u,N}|((0,LN^{n-1}))~\leq~2V\left({N\over L}\right)^{n-1}\,.
 \eeq
\quad\\
 {\bf 3.} Let's define
 \bel{LbN}
 L_N~:=~L\cdot N^{n-1},\qquad\qquad \beta_N~:=~2V\left({N\over L}\right)^{n-1}\,.
 \eeq
We introduce the set
 \begin{multline*}
 \tilde{\mathcal{F}}_N~=~\Big\{f:\left[0, L_N\right]\to [-M,M]~\big|~|Df|((0,L_N))~\leq~\beta_N~\mathrm{and}
 \\
 f(x)=f\left({i\cdot L\over N}\right)\quad\forall x\in \left[{i\cdot L\over N}, {(i+1)\cdot L\over N}\right) \Big\}\,.
 \end{multline*}
From (\ref{TV2}), one has 
 \[
 f_{u,N}~\in~\tilde{\mathcal{F}}_N\qquad\forall u\in \mathcal{F}_{[L,M,V]}\,.
 \]
On the other hand, recalling that 
 \[
\mathcal{B}_{[L_N,2M,\beta_N]}~=~\left\{f\in {\bf L}^1([0,L_N],[0,2M])~\Big|~ |Df|((0,L_N))\leq \beta_N\right\}\,,
 \]
we have 
\[
 \tilde{\mathcal{F}}_N~\subset~\mathcal{B}_{[L_N,2M,\beta_N]}-M\,.
 \]
From Lemma \ref{1-D-BV},  for every $0<\ve'<{L_N\cdot(\beta_N+2M)\over 6}$, it holds
 \[
\mathcal{H}_{\ve'}\left(\mathcal{B}_{[L_N,2M,\beta_N]}~\Big|~{\bf L}^1([0,L_N])\right)~\leq~8\cdot \left\lfloor{L_N(\beta_N+2M)\over \ve'}\right\rfloor\,,
\]
and it yields
\[
\mathcal{H}_{\ve'}\left(\tilde{F}_{N}~\Big|~{\bf L}^1([0,L_N])\right)~\leq~8\cdot \left\lfloor{L_N(\beta_N+2M)\over \ve'}\right\rfloor\,.
\]
By the definition \ref{def1}, there exists a set of $\Gamma_{N,\ve'}=\ds 2^{8\cdot \left\lfloor{L_N(\beta_N+2M)\over \ve'}\right\rfloor}$ functions in $\tilde{F}_N$,
\[
\mathcal{G}_{N,\ve'}~=~\left\{g_{1},g_2,\dots, g_{\Gamma_{N,\ve'}}\right\}~\subset~\tilde{F}_N\,,
\]
such that 
\[
\tilde{F}_{N}~\subset~\ds\bigcup_{i=1}^{\Gamma_{N,\ve'}}~B(g_i,2\ve')\,.
\]
So for every $u \in \mathcal{F}_{[L,M,V]}$, for its corresponding $f_{u,N}, ~\exists~g_{i_u} \in \mathcal{G}_{N,\ve'}$ such that
$$
\|f_{u,N} - g_{i_u}\|_{{\bf L}^1([0,L_N])}~\leq~2\ve'\,.
$$
Let  $\mathcal{U}_{N,\ve'}$ be a set of $\Gamma_{N,\ve'}$ functions $u_j^{\dagger}:[0,L]^N\to[-M,M]$ defined as follows
\begin{eqnarray*}
u^{\dagger}_{j}~=~\left\{\bega{rl}
&\ds 0\qquad\qquad\qquad\qquad~~\mathrm{if}\qquad x\in \bigcup_{\iota\in \{1,2,...,N\}^n}\partial\square_{\iota}\,,
\\[3mm]
&\ds g_{j}\left({{i \cdot L}\over{N}}\right)\quad\qquad\quad~\mathrm{if}\qquad x\in\mathrm{int}\left(\square_{\kappa^i}\right), i\in \{1,2,\dots, N^n\}\,.
\enda\right. 
\end{eqnarray*}
Then corresponding to every $u \in \mathcal{F}_{[L,M,V]}$, there exists $u^{\dagger}_{i_u}\in \mathcal{U}_{N,\ve'}$ for some $i_u\in\{1,2,\dots, \Gamma_{N,\ve}\}$ such that 
\begin{multline*}
\big\|\tilde{u}-u^{\dagger}_{i_u}\big\|_{{\bf L}^1([0,L]^n)}~=~\sum_{i=1}^{N^n}~\left|u_{\kappa^i}- g_{i_u}\left({{i \cdot L}\over{N}}\right)\right|\cdot\mathrm{Vol}\left(\square_{\kappa^i}\right)\\
\qquad~=~\sum_{i=1}^{N^n}~\left| f_{u,N}\left({{i \cdot L}\over{N}}\right)- g_{i_u}\left({{i \cdot L}\over{N}}\right)\right|\cdot {L\over N}\cdot {L^{n-1}\over N^{n-1}}\\
~=~{L^{n-1}\over N^{n-1}}\cdot \|f_{u,N} - g_{i_u}\|_{{\bf L}^1([0,L_N])}~\leq~2\ve'\cdot {L^{n-1}\over N^{n-1}}\,.
\end{multline*}
Combining with (\ref{L-est1}), we  obtain 
\bel{ess2}
\big\|u-u^{\dagger}_{i_u}\big\|_{{\bf L}^1([0,L]^n)}~\leq~\big\|u-\tilde{u}\big\|_{{\bf L}^1([0,L]^n)}+\big\|\tilde{u}-u_{g_{i_u}}\big\|_{{\bf L}^1([0,L]^n)}~\leq~2\ve'\cdot {L^{n-1}\over N^{n-1}}+{L\sqrt{n}\over N}\cdot V\,.
\eeq

{\bf 4.} For any  $\ve>0$, we choose
\bel{cho1}
N~=~\left\lfloor{2\sqrt{n}LV\over \ve}\right\rfloor+1\qquad\mathrm{and}\qquad \ve'~=~{N^{n-1}\cdot \ve\over 4L^{n-1}}~~
\eeq
such that 
\[
\big\|u-u^{\dagger}\big\|_{{\bf L}^1([0,L]^n)}~\leq~2\ve'\cdot {L^{n-1}\over N^{n-1}}+{L\sqrt{n}\over N}\cdot V~\leq~{\ve\over 2}+{\ve\over 2}~=~\ve
\]
for all $u\in \mathcal{F}_{[L,M,V]}$ and for some $u^{\dagger}\in \mathcal{U}_{N,\ve'}$. From the previous step, it holds
\[
\mathcal{F}_{[L,M,V]}~\subseteq~\bigcup_{u^{\dagger}\in \mathcal{U}_{N,\ve'}}~\overline{B}(u^{\dagger},\ve)
\]
provided that 
\bel{cond1}
\ve'~=~{N^{n-1}\ve\over 4L^{n-1}}~\leq~{L_N\cdot(\beta_N+2M)\over 6}~=~{N^{n-1}(VN^{n-1}+ML^{n-1})\over 3 L^{n-2}}\,.
\eeq
This condition is equivalent to 
\[
\ve~\leq~{4\over 3}\cdot \left(LVN^{n-1}+ML^{n}\right)
\]
From (\ref{cho1}), one has that the condition (\ref{cond1}) holds if 
\bel{cd2}
\ve~~\leq~{4\over 3}\cdot \left({2^{n-1}n^{{n-1\over 2}}L^{n}V^n\over \ve^{n-1}}+ML^n\right)\,.
\eeq
Assume that  $0<\ve<{2ML^n\over 3}+ n^{{n-1\over 2n}} LV$, we claim that (\ref{cond1}) holds. Indeed, if ${2ML^n\over 3}> n^{{n-1\over 2n}} LV$ then 
\[
\ve~<~{2ML^n\over 3}+n^{{n-1\over 2n}} LV~\leq~\ds{4ML^n\over 3}
\]
and it yields (\ref{cd2}). Otherwise, we have that  $\ve<{2ML^n\over 3}+n^{{n-1\over 2n}} LV\leq 2n^{{n-1\over 2n}} LV$. Thus
\begin{eqnarray*}
{4\over 3}\cdot \left({2^{n-1}n^{{n-1\over 2}}L^{n}V^n\over \ve^{n-1}}+ML^n\right)&\geq&{4\over 3}\cdot {2^{n-1}n^{{n-1\over 2}}L^{n}V^n\over 2^{n-1}n^{{(n-1)^2\over 2n}}L^{n-1}V^{n-1}}+{4\over 3}ML^n\\[4mm]
&=&{4\over 3}\cdot n^{{n-1}\over 2n}LV+{4\over 3}ML^n\,.
\end{eqnarray*}
and this implies (\ref{cd2}). 

To complete the proof, recalling (\ref{LbN}) and (\ref{cho1}), we estimate 
\begin{eqnarray*}
\mathrm{card}(\mathcal{U}_{N,\ve'})&=&\Gamma_{N,\ve'}~=~2^{8\cdot \left\lfloor{L_N(\beta_N+2M)\over \ve'}\right\rfloor}~=~\ds 2^{8\cdot \left\lfloor{8\over \ve}\cdot \left(LVN^{n-1}+ML^n\right)\right\rfloor}\\[4mm]
&\leq&\ds 2^{{64\over \ve}\cdot\left(LV \left(\left\lfloor{2\sqrt{n}LV\over \ve}\right\rfloor+1\right)^{n-1} +ML^n\right)}\,.
\end{eqnarray*}
Therefore,
\begin{eqnarray*}
\mathcal{H}_{\ve}\left(\mathcal{F}_{[L,M,V]}~\Big|~{\bf L}^1([0,L]^n)\right)&\leq&{64\over \ve}\cdot\left(LV \left(\left\lfloor{2\sqrt{n}LV\over \ve}\right\rfloor+1\right)^{n-1} +ML^n\right)\\[4mm]
&\leq& \ds {64\over\ve}\cdot \left(LV\left({2^{2n-3}n^{{n-1\over 2}}L^{n-1}V^{n-1}\over \ve^{n-1}}+2^{n-2}\right)+ML^n\right)\\[4mm]
&=&\ds{2^{2n+3}n^{n-1\over 2}L^nV^n\over \ve^n}+{2^{n+4}LV+ML^n\over \ve}\,.
\end{eqnarray*}
In particular, if $0<\ve<{ML^n\over 8}$ then 
\[
\mathcal{H}_{\ve}\left(\mathcal{F}_{[L,M,V]}~\Big|~{\bf L}^1([0,L]^n)\right)~\leq~\ds\left[2^{2n+3}n^{n-1\over 2}L^nV^n+\left(2^{n+4}LV+ML^n\right)\cdot\left({ML^n\over 8}\right)^{n-1} \right]\cdot {1\over \ve^n}
\]

and it yields the right hand side of (\ref{m-est}).
\quad\\
\quad\\
{\it (Lower estimate)} We are now going to prove the lower estimate of $\mathcal{H}_{\ve}\left(\mathcal{F}_{[L,M,V]}~\Big|~{\bf L}^1([0,L]^n)\right)$. 

{\bf 1.} Again given any  $N\in\mathbb{N}$, we divide the square $[0,L]^n$ into $N^n$ small squares $\square_{\iota}$ for $\iota=(\iota_1,\iota_2,...,\iota_n)\in \{0,1,...,N-1\}^n$ such that 
\[
\square_{\iota}~=~{\iota L\over N}+\Bigg(\left[0, {L\over N}\right]\times \left[0, {L\over N}\right]\times...\times  \left[0, {L\over N}\right]\Bigg)\qquad\mathrm{and}\qquad \bigcup_{\iota\in \{0,1,2,...,N-1\}^n}~\square_{\iota}~=~[0,L]^n\,.
\]
Consider the set of $N^{n}$-tuples
\[
\Delta_{N}~=~\left\{\ds\delta=(\delta_{\iota})_{\iota\in\{0,1,\dots,N-1\}^n}~\Big|~\delta_{\iota}\in \{0,1\}\right\}\,.
\]
Given any  $h>0$, for any $\delta\in \Delta_N$, define the function $u_{\delta}:[0,L]^n\to \{0,h\}$ such that 
\[
u_{\delta}(x)~=~\sum_{\iota\in \{0,1,\dots,N-1\}^n}h\delta_{\iota}\cdot \chi_{\mathrm{int}(\square_{\iota})}(x)\qquad\forall x\in [0,L]^n\,.
\]
One has $u_{\delta}\in BV((0,L)^n)$ and 
\[
\left|Du_{\delta}\right|((0,L)^n)~\leq~\sum_{\iota\in\{0,1,\dots,N-1\}^n}|Du_{\delta}|(\square_{\iota})~\leq~2^{n-1}\left({L\over N}\right)^{n-1}N^nh~=~(2L)^{n-1}Nh\,.
\]
Assuming that  
\bel{condh}
0~<~h~\leq~\min\left\{M~,~{V\over 2^{n-1}L^{n-1} N}\right\}\,,
\eeq
we have 
\[
\left|Du_{\delta}\right|((0,L)^n)~\leq~(2L)^{n-1}N\cdot {V\over 2^{n-1}L^{n-1} N}~=~V\qquad\forall \delta\in \Delta_N\,,
\]
and this implies 
\[
\mathcal{G}_{h,N}~:=~\left\{u_{\delta}~|~\delta\in\Delta_N\right\}~\subset \mathcal{F}_{[L,M,V]}\qquad\forall N\in\mathbb{N}\,.
\]
Hence,
\bel{comp}
\mathcal{N}_{\ve}\left(\mathcal{F}_{[L,M,V]}~\Big|~{\bf L}^1([0,L]^n)\right)~\geq~\mathcal{N}_{\ve}\left(\mathcal{G}_{h,N}~\Big|~{\bf L}^1([0,L]^n)\right)\qquad\forall \ve>0\,.
\eeq
Towards an estimate of the covering number $\mathcal{N}_{\ve}\left(\mathcal{G}_{h,N}~\Big|~{\bf L}^1([0,L]^n)\right)$, for a fixed $\tilde{\delta}\in\Delta_N$, we can define 
\bel{II}
\mathcal{I}_{\tilde{\delta},N}(2\ve)~=~\left\{\delta\in\Delta_N~\Big|~\|u_{\delta}-u_{\tilde{\delta}}\|_{{\bf L}^1([0,L]^n)}~\leq~2\ve\right\}\qquad\mathrm{and}\qquad C_{N}(2\ve)~=~\mathrm{Card}(\mathcal{I}_{\tilde{\delta},N}(2\ve))
\eeq
since  the cardinality of the set $\mathcal{I}_{\tilde{\delta},N}(\ve)$ is is independent of the choice $\tilde{\delta}\in\Delta_N$. Observe that an $\ve$-cover in ${\bf L}^1$ of  $\mathcal{G}_{h,N}$ contains at most $C_N(2\ve)$ elements. Since $\mathrm{Card}(\mathcal{G}_{h,N})=\mathrm{Card}(\Delta_N)=2^{N^n}$, it holds
\bel{lbb}
\mathcal{N}_{\ve}\left(\mathcal{G}_{h,N}~\Big|~{\bf L}^1([0,L]^n)\right)~\geq~{2^{N^n}\over C_N(2\ve) }\,.
\eeq
{\bf 2.} We now provide an upper bound on $C_N(2\ve)$. For any given pair $\delta,\tilde{\delta}\in\Delta_N$, one has 
\[
\|u_{\delta}-u_{\bar{\delta}}\|_{{\bf L}^1([0,L]^n)}~=~\sum_{\iota\in\{0,1,\dots,N\}^n}\|u_{\delta}-u_{\bar{\delta}}\|_{{\bf L}^1(\square_{\iota})}~=~ d(\delta,\tilde{\delta})\cdot {hL^n\over N^n}\,.
\]
where 
$$
d(\delta,\tilde{\delta})~:=~\mathrm{Card}\left(\{\iota\in\{0,1,\dots,N-1\}^n~|~\delta_{\iota}\neq \tilde{\delta}_{\iota}\}\right)\,.
$$
From (\ref{II}), we obtain 
\[
\mathcal{I}_{\tilde{\delta},N}(2\ve)~=~\left\{\delta\in\Delta_N~\Big|~d(\delta,\tilde{\delta})~\leq~{2\ve N^n\over hL^n}\right\}\,,
\]
and it yields
\[
C_N(2\ve)~=~\mathrm{Card}\left(\mathcal{I}_{\tilde{\delta},N}(2\ve)\right)~\leq~\sum\limits_{r = 0}^{\left\lfloor{2\ve N^n\over hL^n}\right\rfloor} {N^n \choose r}\,.
\]
To estimate the last term in the above inequality, let's consider $N^n$ independent random variables with uniform Bernoulli distribution  $X_1,X_2,\dots, X_{N^n}$
\[
\mathbb{P}(X_i=1)~=~\mathbb{P}(X_i=0)~=~{1\over 2}\qquad\forall i\in \{1,2,\dots, N^n\}\,.
\]
Set $S_{N^n}:= X_1+X_2+\dots+X_{N^n}$. Observe that for any $k\leq N^n$, we have 
\[
\sum_{r=1}^{k}~ {N^n \choose r}~=~2^{N^n}\cdot \mathbb{P}\left(S_{N^n}\leq k\right)\,.
\]
Thanks to Hoeffding's inequality \cite[Theorem]{Hoeffding}, for all $\mu\leq {N^n\over 2}$, one has
\[
\mathbb{P}\left(S_{N^n}\leq \mathbb{E}[S_{N^n}]-\mu\right)~=~\mathbb{P}\left(S_{N^n}\leq {N^n\over 2}-\mu\right)~\leq~\exp\left(-{2\mu^2\over N^n}\right)
\]
where $ \mathbb{E}[S_{N^n}]$ is the expectation of $S_{N^n}$.
Hence, for every $0<\ve\leq{hL^n\over 8}$ such that ${2\ve N^n\over hL^n}\leq {N^n\over 2}$ and ${4\ve\over hL^n}\leq {1\over 2}$, it holds
\begin{eqnarray*}
C_N(2\ve)&\leq&\sum\limits_{r = 0}^{\left\lfloor{2\ve N^n\over hL^n}\right\rfloor} {N^n \choose r}~=~ 2^{N^n}\cdot\mathbb{P}\left(S_{N^n}\leq \left\lfloor{2\ve N^n\over hL^n}\right\rfloor\right)
\\[4mm]
&\leq&2^{N^n}\cdot \exp\left(-{2\left({N^n\over 2}-\left\lfloor{2\ve N^n\over hL^n}\right\rfloor\right)^2\over N^n }\right)~\leq~2^{N^n}\cdot \exp\left(-{\left(N^n-{4\ve N^n\over hL^n}\right)^2\over 2N^n }\right)\\[4mm]
&=&2^{N^n}\cdot \exp\left(-N^n\cdot {\left(1-{4\ve\over hL^n}\right)^2\over 2}\right)~\leq~2^{N^n}\cdot e^{-N^n/8}\,.
\end{eqnarray*}
From (\ref{lbb}) and (\ref{condh}), the following holds
\begin{eqnarray*}
\mathcal{N}_{\ve}\left(\mathcal{G}_{h,N}~\Big|~{\bf L}^1([0,L]^n)\right)&\geq&{2^{N^n}\over C_N(2\ve)}~\geq~e^{{N^n\over 8}}
\end{eqnarray*}
provided that 
\bel{condl}
0~<~h~\leq~\min\left\{M~,~{V\over 2^{n-1}L^{n-1} N}\right\}\qquad\mathrm{and}\qquad 0<\ve~\leq~{hL^n\over 8}\,.
\eeq
Therefore, for every $0<\ve<{ML^n\over 8}$, by choosing 
\[
 h~=~\min\left\{M~,~{V\over 2^{n-1}L^{n-1} N}\right\}\qquad\mathrm{and}\qquad N~\doteq~\left\lfloor{VL\over 2^{n+2}\ve}\right\rfloor
\]
such that (\ref{condl}) holds, we obtain that
\[
\mathcal{N}_{\ve}\left(\mathcal{G}_{h,N}~\Big|~{\bf L}^1([0,L]^n)\right)~\geq~\exp\left({1\over 8}\cdot \left\lfloor{VL\over 2^{n+2}\ve}\right\rfloor^n\right)\,.
\]
Recalling (\ref{comp}), we have 
\[
\mathcal{N}_{\ve}\left(\mathcal{F}_{[L,M,V]}~\Big|~{\bf L}^1([0,L]^n)\right)~\geq~\exp\left({1\over 8}\cdot \left\lfloor{VL\over 2^{n+2}\ve}\right\rfloor^n\right)
\]
and this implies the first inequality in (\ref{m-est}).
\endproof

\v
{\bf Acknowledgments.} K.T. Nguyen is partially supported by a grant from the Simons Foundation/SFARI (521811, NTK).

\end{document}